\documentclass{ifacconf}

\usepackage{natbib}            
\usepackage{graphicx} 
\usepackage[utf8]{inputenc}
\usepackage[T1]{fontenc}
\usepackage[retainorgcmds]{IEEEtrantools}
\usepackage{amssymb,amsmath}
\usepackage{url}
\usepackage{siunitx}
\usepackage{xcolor}

\graphicspath{{Figures/}}

\begin{document}
	
\begin{frontmatter}
	
\title{\Large Continuous globally exponential rigid-body angular velocity estimation from attitude measurements}

\author[First]{Ioannis Sarras}
		
\address[First]{Office National d'\'Etudes et de Recherches A\'erospatiales (ONERA)--The French Aerospace Lab, F-91123 Palaiseau, France (e-mail: ioannis.sarras@onera.fr).}
		
				
\begin{abstract}
A continuous solution is proposed to the problem of uniform global exponential estimation of the angular velocity for rigid bodies by means of direct attitude measurements. The proposed observer is designed on  $\mathbb{R}^3\times\mathbb{R}^3\times\mathbb{R}^3$, and not on the inherent geometry $\mathtt{SO}(3)\times\mathbb{R}^3$, which allows for a simple structure of the observer, being a copy of the system plus some correction terms, and a straightforward classical Lyapunov stability proof. The performances of the observer are illustrated by means of a simple simulation scenario.
\end{abstract}
		
\end{frontmatter}

\section{Introduction}\label{sec:intro}

The problem of attitude control for rigid bodies has been driving for decades the developments in both theoretical and applied research due to its rich characteristics and major industrial impact. Throughout the years a variety of control algorithms have been proposed that can in one way be classified into velocity-independent and velocity-dependent based on whether they require measurements of angular velocity or not. The first category exploits, the well-known by now, passivity properties of the rigid-body dynamics  to use only orientation information, see e.g. \cite{LizarW1996TAC,Tsiot1998TAC}, while the second category adopts an observer-controller structure with the observer providing an estimate of the angular velocity by exploiting direct attitude measurements, see \cite{AkellT2015JGCD,BerkaAT2016CDC} and references therein for an exhaustive list of references.  

Here the focus is on the second category and more precisely, we consider the problem of designing a \emph{continuous} full-order angular velocity observer that provides a \emph{uniformly global exponential estimate}. This question has been adressed by various researchers since the seminal paper of \cite{Salcu1991TAC} where an observer exploiting the quaternion formulation was employed to generate a (non uniform) globally convergent estimate. Recently, several works \cite{AkellT2015JGCD,BerkaAT2016CDC} have shown anew interest for this problem in order to provide for the first time uniform and global, asymptotic or exponential, solutions along with a (strict) Lyapunov stability proof which in turn can serve to establish among others stability of observer-controller interconnections and robustness margins with respect to disturbances.

Compared to the recent work by \cite{BerkaAT2016CDC}, the proposed observer of this work is continuous and thus, does not require any hysterisis-based mechanism to switch between different configurations. Furthermore, the stability analysis presented here hinges upon classic Lyapunov arguments and does not require advanced tools from hybrid systems theory in order to conclude the uniformly global exponential stability claim. Compared to the other recent work of \cite{AkellT2015JGCD}, we prove uniformly global exponential stability of zero estimation error instead of just uniformly global asymptotic stability. In addition, through an appropriate error definition a strict Lyapunov function is directly provided and hence, no strictification procedure is required as in \cite{AkellT2015JGCD} that furthermore provides a partial strict Lyapunov function.

The proposed observer does not respect the inherent geometry of the problem, that is does not live in $\mathtt{SO}(3)\times \mathbb{R}^3$, but is rather designed on  $\mathbb{R}^3\times\mathbb{R}^3\times \mathbb{R}^3$. The adopted parameterization, inspired by the developments in \cite{BatisSO2012AUT,BatisSO2012CEP,BatisSO2014CDC,MartiS2016AUT,MartiS2016CDC,MartiS2017IFAC}, leads to a simple structure of the observer, being a copy of the system plus some nonlinear correction terms, and a rather direct stability proof based on a strict Lyapunov function.

The design of the observer hinges upon the invariant-manifold based methodology introduced in \cite{KaragSA2009AUT}, see also \cite{InIbook}. 
Apart from the applications in the above references, this methodology has been recently applied to succesfully solve a variety of problems such as global velocity estimation for mechanical systems \cite{AstolOV2013AUT}, \cite{VenkaOSV2010TAC}, global exponential position-feedback tracking for fully-actuated mechanical systems \cite{RomerOS2015TAC}, global asymptotic position-feedback synchronization for tele-operation systems \cite{SarraNKB2016IJRNC}, as well as for semi-global reduced attitude estimation for quadrotors \cite{MartiS2016CDC} and global exponential velocity-aided full attitude  estimation \cite{MartiS2016AUT} for rigid bodies.

The presentation of the paper is organized as follows. Section \ref{sec:model} presents the model on $\mathtt{SO}(3)\times \mathbb{R}^3$, for the attitude kinematics and the angular dynamics, and the representation on $\mathbb{S}^2\times\mathbb{S}^2\times \mathbb{R}^3$ that will be adopted for the observer design. Section \ref{sec:observer} describes the main contribution which is the design of a full-order observer for the angular velocity on $\mathbb{R}^3\times\mathbb{R}^3\times \mathbb{R}^3$ and the corresponding stability proof. Then follows Section \ref{sec:simu} that illustrates the performances of the proposed estimator in a simulation scenario. The paper is concluded with some conclusions and perspectives.  
\section{Model and Measurements} \label{sec:model}

We consider a moving rigid body subjected to the angular velocity~$\omega\in\mathbb{R}^3$ (in body axes). Its orientation (from inertial to body axes) matrix $R\in\mathtt{SO}(3)$ is related to~$\omega$ through the orientation kinematics as

\begin{IEEEeqnarray}{rCl}
	\dot R &=& R\omega_\times\label{eq:R},
\end{IEEEeqnarray} 

where $\omega_\times$ is the skew-symmetric matrix defined by $\omega_\times x:=\omega\times x$ whatever the vector~$x$. The attitude kinematics is complemented by the attitude dynamics that is given by 

\begin{IEEEeqnarray}{rCl}
	J\dot\omega &=& (J\omega)\times\omega + \tau,
\end{IEEEeqnarray}

with $J\in\mathbb{R}^{3\times 3}$ the rigid body inertia matrix and $\tau\in\mathbb{R}^3$ the input torque. 

Inspired by the recent works on attitude estimation \cite{BatisSO2012AUT,BatisSO2012CEP,BatisSO2014CDC,MartiS2016AUT,MartiS2016CDC,MartiS2017IFAC}, the following parametrization of the attitude matrix $R$ is adopted

\begin{IEEEeqnarray}{rCl}
R	=\left[\begin{array}{c}
r_1^T\\
r_2^T\\
r_3^T
\end{array}\right],
\end{IEEEeqnarray}

with $r_i\in\mathbb{S}^2,~i=\{1,2,3\}$. Then, and using the kinematics \eqref{eq:R}, one can express the $r_i$-dynamics as

\begin{IEEEeqnarray}{rCl}
	\dot r_1 &=& r_1\times\omega\\
	\dot r_2 &=& r_2\times\omega\\
	\dot r_3 &=& r_3\times\omega\\
	J\dot\omega &=& (J\omega)\times\omega + \tau\\
	y&=&\left[\begin{array}{c}
		r_1\\
		r_2\\
		r_3
	\end{array}\right].
\end{IEEEeqnarray}

Since the rotation matrix $R$ is orthogonal with orthonormal rows $r_i$, it is evident that one only requires the knowledge of two of the rows of $R$ to reconstruct $\omega$, since the third row can be obtained by the cross-product of the other two. As such, and without loss of generality, it will be considered henceforth that only the vectors $r_1,r_2$ are measured and that these will be used for the observer design. 

Then, the considered model is reduced to 

\begin{IEEEeqnarray}{rCl}
	\dot r_1 &=& r_1\times\omega, \quad |r_1|=1\\ \label{eq:r1}
	\dot r_2 &=& r_2\times\omega,, \quad |r_2|=1\\
	J\dot\omega &=& (J\omega)\times\omega + \tau\label{eq:w}\\
	y&=&\left[\begin{array}{c}
		r_1\\
		r_2
	\end{array}\right].
\end{IEEEeqnarray}

Having presented the underlying model, whose underlying geometry is in $\mathbb{S}^2\times\mathbb{S}^2\times \mathbb{R}^3$, we can now state the working assumptions on which the observer design will be based.

\begin{hypo}[A.1]
	The inertia matrix $J$ is symmetric, positive definite and uniformly bounded, i.e. $J_mI\leq J\leq J_MI$ with known bounds $J_m,J_M>0$.\label{as:1}
\end{hypo}
\begin{hypo}[A.2]
	The angular velocity $\omega$ is uniformly bounded, i.e. $|\omega|\leq \omega_M$ with a known bound $\omega_M>0$.\label{as:2}
\end{hypo}

Furthermore, a result that is crucial for the stability analysis is reminded.
\begin{fact}[\cite{BenziBCT2015IJRNC,TayebRB2013TAC}]\label{fact:1}
	~\\
	Assume that two vectors $r_1,r_2\in\mathbb{R}^3$ are linearly independent. Then, the matrix 
	
	\begin{IEEEeqnarray}{rCl}
	M(r_1,r_2)&:=& -k_1 r_{1\times}^2 - k_2 r_{2\times}^2,
	\end{IEEEeqnarray}

is symmetric and positive definite, i.e. $M=M^T\succ 0$. Furthermore, the positive scalars $k_1$, $k_2$ can freely increase the (simple) eigenvalues of $M(r_1,r_2)$.
\end{fact}

This section is concluded by stating the problem to be solved.

\emph{Problem statement.} Design a continuous observer based on the model \eqref{eq:r1}-\eqref{eq:w}, measurements $r_1,~r_2$ and under assumptions A.~\ref{as:1}-A.~\ref{as:2}, that provides a \emph{uniformly globally exponentially convergent} angular velocity estimate. 

\section{Observer design}\label{sec:observer}
The main result of this paper can now be stated.

\begin{thm}
Consider the dynamics~\eqref{eq:r1}-\eqref{eq:w} and the working assumptions A.\ref{as:1}-A.\ref{as:2}. Then the dynamical system
 
	\begin{IEEEeqnarray}{rCl}
		\dot \xi &=& -k_1J^{-1}(r_{1\times}\hat{r}_{1\times}-\hat{r}_{1\times}r_{1\times})\hat{\omega}\nonumber\\
		&-&k_2J^{-1}(r_{2\times}\hat{r}_{2\times}-\hat{r}_{2\times}r_{2\times})\hat{\omega}\nonumber\\
		&+&l_1(\hat{\omega})k_1J^{-1}r_1\times(\hat{r}_1-r_1)+l_2(\hat{\omega})k_2J^{-1}r_2\times(\hat{r}_2-r_2)\nonumber\\
		&+&J^{-1}((J\hat\omega)\times\hat\omega +\tau)\\
		\hat{\omega}&:=&\xi+J^{-1}k_1\hat{r}_1\times  r_1+J^{-1}k_2\hat{r}_2\times r_2\\
		\dot{\hat{r}}_1 &=& \hat{r}_1\times\hat\omega -l_1(\hat{\omega})(\hat{r}_1-r_1)\\
		\dot{\hat{r}}_2 &=& \hat{r}_2\times\hat\omega -l_2(\hat{\omega})(\hat{r}_2-r_2)
	\end{IEEEeqnarray}
is a uniformly global exponential observer under the following conditions on the gains

\begin{IEEEeqnarray}{rCl}
	-k_1 r_{1\times}^2 &-& k_2 r_{2\times}^2\succ(\omega_M+\frac{3}{2})I\label{gain1}\\
	l_1&>&(|\hat{\omega}|+\omega_M)(1+\frac{k_1}{J_m}(1+\frac{2k_1}{J_m})\\
	l_2&>&(|\hat{\omega}|+\omega_M)(1+\frac{k_2}{J_m}(1+\frac{2k_2}{J_m})\label{gain3}.
\end{IEEEeqnarray}	

\end{thm}

\begin{pf}
Let us define the estimation error

\begin{IEEEeqnarray*}{rCl}
	z &=& \xi +\beta(y,\hat y) - \omega=:\hat{\omega}-\omega,
\end{IEEEeqnarray*}
with the mapping $\beta(y,\hat y)$ defined as

\begin{IEEEeqnarray*}{rCl}
	\beta(y,\hat y)&:=&J^{-1}k_1\hat{r}_1\times r_1+J^{-1}k_2\hat{r}_2\times r_2.
\end{IEEEeqnarray*}
The error dynamics then can be expressed as

\begin{IEEEeqnarray*}{rCl}
	\dot z &=& \dot\xi +\partial_y\beta\dot y +\partial_{\hat y}\beta \dot{\hat y}- \dot\omega\\
	&=& \dot\xi +\partial_y\beta\left[\begin{array}{c}
		r_{1\times}\\
		r_{2\times}
	\end{array}\right]\omega - J^{-1}((J\omega)\times\omega+\tau)\\
	&=&-\partial_y\beta\left[\begin{array}{c}
		r_{1\times}\\
		r_{2\times}
	\end{array}\right]z + J^{-1}((J\hat\omega)\times\hat\omega - (J\omega)\times\omega)\\
	&=& -\partial_y\beta\left[\begin{array}{c}
		r_{1\times}\\
		r_{2\times}
	\end{array}\right]z + J^{-1}((Jz)\times z + (J\omega)\times z + (Jz)\times \omega)
\end{IEEEeqnarray*}
with the observer dynamics' general form chosen as

\begin{IEEEeqnarray}{rCl}
	\dot \xi &=& -\partial_{\hat y}\beta\dot{\hat y}-\partial_y\beta\left[\begin{array}{c}
		r_{1\times}\\
		r_{2\times}
	\end{array}\right]\hat\omega+ J^{-1}((J\hat\omega)\times\hat\omega +\tau)
\end{IEEEeqnarray}

Using the definition of the mapping $\beta$ one can straightforwardly obtain the explicit form of the $z$-dynamics that gives

\begin{IEEEeqnarray*}{rCl}
	\dot z &=& J^{-1}(k_1\hat{r}_{1\times}r_{1\times} + k_2\hat{r}_{2\times}r_{2\times})z\\
	 &+& J^{-1}((Jz)\times z + (J\omega)\times z + (Jz)\times \omega)\\
	&=& J^{-1}(k_1r_{1\times}^2 + k_2r_{2\times}^2)z\\
	&+& J^{-1}((Jz)\times z + (J\omega)\times z + (Jz)\times \omega)\\
	&+& J^{-1}k_1(e_{1\times}^2+e_{1\times}r_{1\times}+r_{1\times}e_{1\times})z\\
	&+& J^{-1}k_2(e_{2\times}^2+e_{2\times}r_{2\times}+r_{2\times}e_{2\times})z,
\end{IEEEeqnarray*}
where we defined the filtering errors $e_1:=\hat{r}_1-r_1$ and $e_1:=\hat{r}_2-r_2$.

Let us now consider the following Lyapunov function candidate 
\begin{IEEEeqnarray*}{rCl}
	V_z(z) &=& |Jz|.
\end{IEEEeqnarray*}

Its time derivative along $z-$trajectories gives 

\begin{IEEEeqnarray*}{rCl}
	\dot{V}_z &=&  \frac{z^T(k_1r_{1\times}^2 + k_2r_{2\times}^2)z}{|Jz|} +\frac{z^T(Jz)\times \omega}{|Jz|}\\
	&+&\frac{z^T(k_1(e_{1\times}^2+e_{1\times}r_{1\times}+r_{1\times}e_{1\times})))z}{|Jz|}\\
	&+&\frac{z^T(k_2(e_{2\times}^2+e_{2\times}r_{2\times}+r_{2\times}e_{2\times}))z}{|Jz|}\\
	&\leq&-\frac{\mu}{J_M}|z|+\omega_M|z| +\frac{|z|}{J_m}(k_1(|e_1|^2+2|e_1|)\\
	&+&k_2(|e_2|^2+2|e_2|))\\
	&\leq& -(\frac{\mu}{J_M}-\omega_M)|z|+\frac{k_1}{J_m}(|\hat{\omega}|+\omega_M)|e_1|^2 +\frac{1}{2}|z|\\
	&+&\frac{2k_1^2}{J_m}(|\hat{\omega}|+\omega_M)|e_1|^2+\frac{k_2}{J_m}(|\hat{\omega}|+\omega_M)|e_2|^2 \\
	&+&\frac{1}{2}|z|+\frac{2k_2^2}{J_m}(|\hat{\omega}|+\omega_M)|e_2|^2\\
	&\leq& -(\frac{\mu}{J_M}-\omega_M-1)|z|+\frac{k_1}{J_m}(|\hat{\omega}|+\omega_M)(1+\frac{2k_1}{J_m})|e_1|^2\\
	&+&\frac{k_2}{J_m}(|\hat{\omega}|+\omega_M)(1+\frac{2k_2}{J_m})|e_2|^2,
\end{IEEEeqnarray*}

where we have used from Fact~\ref{fact:1} that 

\begin{IEEEeqnarray*}{rCl}
	-k_1 r_{1\times}^2 - k_2 r_{2\times}^2\succeq \mu I\succ0,
\end{IEEEeqnarray*}
with $\mu$ taking larger values for larger values of the assignable gains $k_1,k_2$.

Now, define the filter dynamics

\begin{IEEEeqnarray}{rCl}
	\dot{\hat{r}}_1 &=& \hat{r}_1\times\hat\omega -l_1(\hat{r}_1-r_1)\\
	\dot{\hat{r}}_2 &=& \hat{r}_2\times\hat\omega -l_2(\hat{r}_2-r_2),
\end{IEEEeqnarray}
that lead to the following error dynamics

\begin{IEEEeqnarray*}{rCl}
	\dot{e}_1 &=& r_1\times z + e_1\times\hat\omega -l_1e_1\\
	\dot{e}_2 &=& r_2\times z + e_2\times\hat\omega -l_2e_2.
\end{IEEEeqnarray*}

To establish the stability analysis for the ($e_1,e_2$)-dynamics the positive definite, radially unbounded, function $	V_e(e_1,e_2)$ is defined as

\begin{IEEEeqnarray*}{rCl}
	V_e(e_1,e_2) &=& \frac{1}{2}(|e_1|^2+|e_2|^2).
\end{IEEEeqnarray*}

Its time-derivative along ($e_1,e_2$)-trajectories yields

\begin{IEEEeqnarray*}{rCl}
	\dot V_e&=& e_1^T(r_1\times z)-l_1|e_1|^2+e_2^T(r_2\times z)-l_2|e_2|^2\\
	&\leq& -l_1|e_1|^2-l_2|e_2|^2 + (|e_1|+|e_2|)|z|\\
	&\leq& -l_1|e_1|^2-l_2|e_2|^2 + \frac{|z|}{2}+(|e_1|^2+|e_2^2|)|z|\\
	&\leq& -(l_1-(|\hat{\omega}|+\omega_M))|e_1|^2-(l_2-(|\hat{\omega}|+\omega_M))|e_2|^2\\
	&+& \frac{|z|}{2}.
\end{IEEEeqnarray*}

Finally, putting all pieces together, we have the composite Lyapunov function

\begin{IEEEeqnarray*}{rCl}
	V(z,e_1,e_2) &=& V_z(z)+V_e(e_1,e_2),
\end{IEEEeqnarray*}
whose time-derivative evaluated along $(z,e_1,e_2)$-trajectories gives

\begin{IEEEeqnarray*}{rCl}
	\dot V&\leq& -(\frac{\mu}{J_M}-\omega_M-\frac{3}{2})|z|\\
	&-&(l_1-(|\hat{\omega}|+\omega_M)(1+\frac{k_1}{J_m}(1+\frac{2k_1}{J_m}))|e_1|^2\\
	&-&(l_2-(|\hat{\omega}|+\omega_M)(1+\frac{k_2}{J_m}(1+\frac{2k_2}{J_m}))|e_2|^2\\
	&\leq& -(\mu-J_M(\omega_M+\frac{3}{2}))\frac{|Jz|}{J_M^2}\\
	&-&(l_1-(|\hat{\omega}|+\omega_M)(1+\frac{k_1}{J_m}(1+\frac{2k_1}{J_m}))|e_1|^2\\
	&-&(l_2-(|\hat{\omega}|+\omega_M)(1+\frac{k_2}{J_m}(1+\frac{2k_2}{J_m}))|e_2|^2,
\end{IEEEeqnarray*}

from which the uniform global exponential stability claim follows under the gain conditions \eqref{gain1}-\eqref{gain3}.

\end{pf}
\section{Simulations}\label{sec:simu}

The excellent behavior of the observer is now illustrated in a simple simulation scenario. With respect to the physical parameters we consider an inertia matrix given as $J={\tt diag}(3,2,1)$. As for the system, the initial conditions on the angular velocity are taken as $\omega_0 = {\tt col}(2,-1,0.6)]$~(rad/sec) and we further consider the upper bound on the angular velocity $\omega_M = 3$~(rad/sec). Concerning the initial orientation we consider that $R(0)=I$ and thus, have that $r_{10} = {\tt col}(1,0,0)$ and
$r_{20} = {\tt col}(0,1,0)$.

The inital state of the observer has been chosen in such a way that the initial velocity estimate $\hat{\omega}_0=0$ by taking
$\xi_0= k_1J^{-1}(\hat{r}_{10}\times r_{10})+k_2J^{-1}(\hat{r}_{20}\times r_{20}).$ In addition, the initial estimates $\hat{r}_{10},\hat{r}_{20}$ are selected as
$\hat{r}_{10} = {\tt col}(0,-1,0)$, $\hat{r}_{10} = {\tt col}(0,0,-1)$ in order to illustrate that they are not designed to live on the unitary sphere.

Finally, the tuning gains are set to $(k_1,k_2,l_{10},l_{20})=(14.5,14.5,1,1)$ and where the gains $k_1,k_2$ have been taken as $k_1 = 1+(\omega_M+3/2)J_M$, 
$k_2 = 1+(\omega_M+3/2)J_M$.

\begin{figure}[ht]
	\includegraphics[width=\columnwidth]{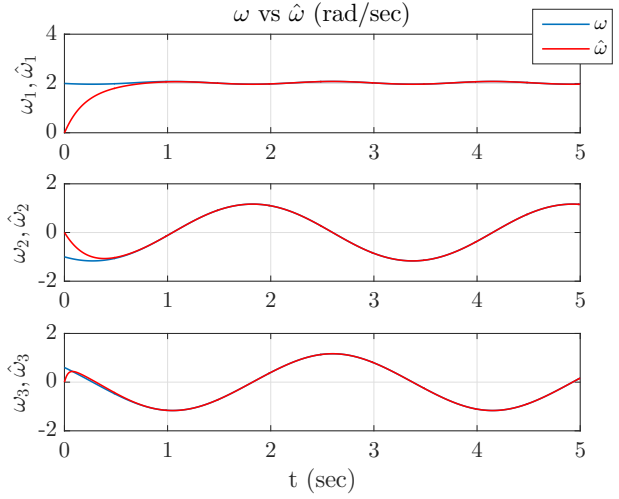}
	\caption{Components of the true $\omega$ (blue) and the estimation $\hat{\omega}$ (red).}
	\label{fig:OmegaError}
\end{figure}

\begin{figure}[ht]
	\includegraphics[width=\columnwidth]{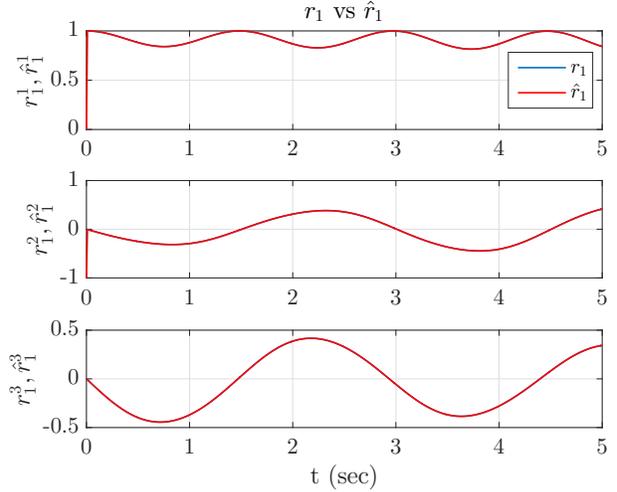}
	\caption{Components of the true $r_1$ (blue) and the filtered $\hat{r}_1$ (red).}
	\label{fig:r1Error}
\end{figure}

The time responses of the angular velocity estimates along with the true angular velocities are depicted in Fig.~\ref{fig:OmegaError}. The transient behavior and convergence is smooth and fast as expected, and despite the initial errors. Larger values of the gains allow for faster convergence.   

\begin{figure}[ht]
	\includegraphics[width=\columnwidth]{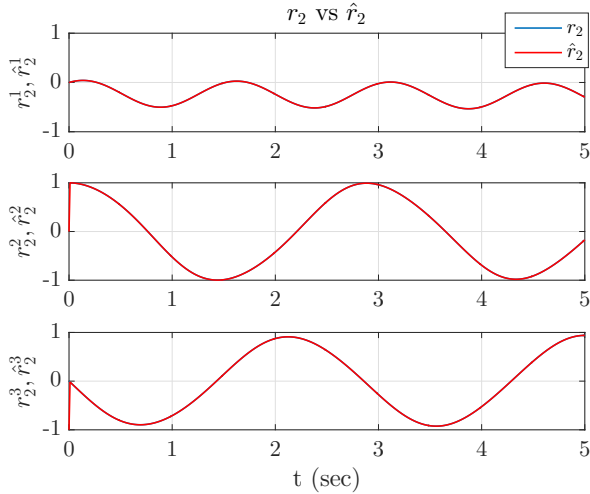}
	\caption{Components of the true $r_2$ (blue) and the filtered $\hat{r}_2$ (red).}
	\label{fig:r2Error}
\end{figure}

In addition, the time evolution of the filtered measurements $\hat{r}_1,\hat{r}_2$ versus the time evolution of the true measurements $r_1,r_2$ is illustrated in Fig.~\ref{fig:r1Error} and  Fig.~\ref{fig:r2Error} respectively. 

\begin{figure}[ht]
	\includegraphics[width=\columnwidth]{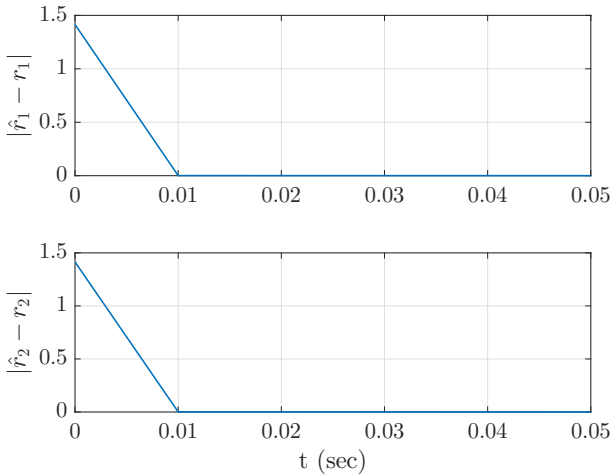}
	\caption{Norm of the filtering errors $|\hat{r}_1-r_1|$ and $|\hat{r}_2-r_2|$.}
	\label{fig:NormFiltErrors}
\end{figure}

Finally, the norm of the filtering errors $\hat{r}_1-r_1$, $\hat{r}_2-r_2$ is presented in Fig.~\ref{fig:NormFiltErrors} in order to show the different time-scale for the convergence of the filters on $r_1,r_2$ with respect to the angular velocity estimate.

\section{Conclusions--Perspectives}
A continuous-time, full-order nonlinear observer has been presented that provides a uniformly globally exponentially convergent estimate of the rigid body's angular velocity from direct attitude measurements. The selected attitude parametrization along with the adopted invariant-manifold--based observer design methodology reveal an appropriate error definition and a simple observer dynamics' structure that allow for a straightforward (strict) Lyapunov-based stability proof. To the best of the author's knowledge this is the first work in the literature where a continuous-time, uniformly global exponential solution is proposed for the problem of interest and moreover, with both a simple observer structure and a simple stability analysis. The performances of the proposed observer are illustrated by means of a realistic simulation scenario.

Given the advantageous characteristics of the presented observer, future promising directions include the extension to the case of a single vector (rotation matrix row) measurement, that will certainly require a persistence-of-excitation assumption, as well as the study of a controller-observer interconnection for global attitude stabilization/tracking.

\end{document}